\documentclass[11pt,twoside]{article}
\usepackage{latexsym,amssymb,amsfonts}
\usepackage{amsmath}
\usepackage[english]{babel}
\usepackage{epsfig}
\begin{document}

\vskip -3.0in

\title{Quantum deformations of fundamental groups of oriented $3$-manifolds}

\author{Uwe Kaiser \\
\begin{small}
\makeatletter
\centerline{e-mail: kaiser@math.boisestate.edu}
\makeatletter
\end{small}
}

\maketitle

\begin{small}
\vskip .1in

\noindent \textsf{ABSTRACT.} We compute two-term skein modules of framed oriented links in oriented $3$-manifolds. They contain the self-writhe and total
linking number invariants of framed oriented links in a universal way.
 The relations in a natural presentation of the skein module are interpreted as monodromies in the space of immersions from circles to the $3$-manifold.

\vskip .05in


\vskip .05in

\noindent \textsf{MATHSUBJECT CLASS.} 57M25, 57M35, 57R42
\end{small}

\section{Introduction}

Let $M$ be a compact oriented $3$-manifold (possibly with boundary). Throughout link will mean framed and oriented link in $M$. 

In [P1] (see also [Ki]) J.\ Przytycki defines the \textit{quantum deformation of $\pi_1(M)$} in the following way:
Let $\mathfrak{L}(M)$ be the set of isotopy classes of links in $M$. Let $\mathcal{S}(M)$ be the quotient of the free $\mathbb{Z}[q^{\pm 1}]=:R$-module with basis $\mathfrak{L}(M)$ by the relations $K_+=q^2K_-$ and $K^{(1)}=qK$. Here $K_{\pm}$ are two  links, which differ only inside a $3$-ball in $M$ by a crossing change. The link $K^{(1)}$ is defined from $K$ by introducing a positive twist into the framing of one of its components. The pictures below show projections onto a plane contained in some oriented $3$-ball in $M$. We assume that framings are induced from the projection plane. 

\vskip .2in

\begin{center}
\epsfig{file=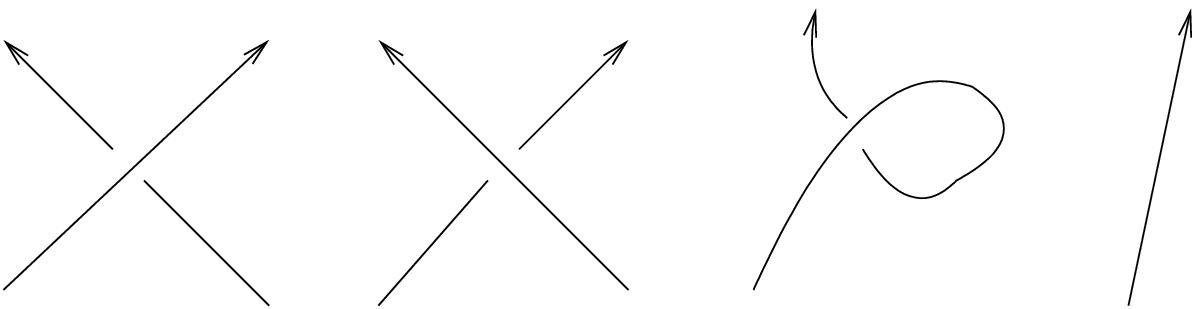,height=1.3in,width=4in}
\end{center}

\qquad \qquad \qquad \qquad \quad $K_+$ \quad \qquad \qquad \qquad $K_-$ \qquad \qquad \qquad $K^{(1)}$ 
\qquad \qquad  $K$ 

\vskip .2in

Let $\hat{\pi}(M)$ be the set of conjugacy classes of the fundamental group $\pi_1(M)$. For each ring $A$ let $SA\hat{\pi}(M)$ be the symmetric algebra of the free module with basis $\hat{\pi }(M)$.

\vskip .1in

\noindent \textbf{Theorem (Przytycki}). \textbf{(a)} \textit{If $M$ has no non-separating $2$-spheres and tori then $\mathcal{S}(M)$ is isomorphic to $SR\hat{\pi}(M)$.}

\noindent \textbf{(b)} \textit{If $M$ contains a non-separating $2$-sphere or torus then $\mathcal{S}(M)$ has a torsion element:}

\begin{enumerate}
\item[\textbf{(i)}] \textit{If $K$ is a link in $M$ with intersection number with some $2$-sphere equal to $k\neq 0$ then \newline  $(q^{2k}-1)K=0$ in $\mathcal{S}(M)$.}
\item[\textbf{(ii)}] \textit{Let $K'$ be a link in $M$ with intersection number with some torus equal to $k\neq 0$. Let $K$ be a link obtained from $K'$ by adding to $K'$ a non-contractible curve on the torus. Then $(q^{2k}-1)K=0$ in 
$\mathcal{S}(M)$.}
\end{enumerate}

The results \textbf{(b)} are proven in ([P1], 1.3)  by explicit construction and \textbf{(a)} has been announced there. We will give a proof of Przytycki's theorem in section 5. 
 A similar and related skein module based on homology of framed links has been computed in [P3]. This module has been described as the \textit{quantum deformation of the first homology group of $M$} in the same way as $\mathcal{S}(M)$ can be considered as the \textit{quantum deformation of the fundamental group of $M$.}
Note that the ring homomorphism $R\rightarrow 1$ mapping $q$ to $1$ induces the $(R,\mathbb{Z})$-homomorphism
$$\mathcal{S}(M)\rightarrow S\mathbb{Z}\hat{\pi}(M).$$

In this paper we give a complete description of $\mathcal{S}(M)$, which is based on ideas from Vassiliev theory, see also [KL],[K1],[V1].
 This settles Problem I.\ 92, Part I from Kirby's problem list. We will reduce the understanding of relations in a natural presentation of $\mathcal{S}(M)$ to certain problems in the theory of \textit{singular} $2$-spheres and tori in $M$.
The relation with Przytycki's result is based on D.\ Gabai's generalization of the sphere theorem relating singular with embedded surface theory [G].

\vskip .1in

For $r\geq 0$ let $\Lambda_r(M):=map(\amalg_rS^1,M)/\Sigma_r$, where $map$ denotes the set of smooth maps and the permutation group $\Sigma_r$ acts on domains in the obvious way. Note that $\Lambda_1(M)=:\mathbb{L}(M)$ is the space of smooth free loops in $M$ and $\Lambda_0(M)$ is a point by definition.
 Let $\Lambda(M):=\amalg_{r\geq 0}\Lambda_r(M)$ be the \textit{generalized free loop space of $M$}. Let 
$I(M)\subset \Lambda(M)$ be the subspace of immersed maps. 

\vskip .1in

We will prove the following two results:

\vskip .1in

\noindent \textbf{Theorem 1.} \textit{There exists a system of local coefficients $\mathfrak{Z}$ in $I(M)$, with $\mathfrak{Z}_x\cong R$ for $x\in I(M)$, such that 
$$H_0(I(M),\mathfrak{Z})\cong \mathcal{S}(M).$$}

\vskip -.1in

 Note that $\pi_0(\mathbb{L}(M))$ is in 1-1 correspondence with  
$\hat{\pi}(M)$ and
 $\pi_0(\Lambda(M)) \cong \pi_0(I(M))$ is in 1-1 correspondence with the set 
$\mathfrak{b}(M)$ of unordered sequences in $\hat{\pi }(M)$. Also $\mathfrak{b}(M)$ is the natural basis of the module $SA\hat{\pi }(M)$, which is isomorphic to the homology module $H_0(\Lambda (M),A)$ ($A$ any ring).

\vskip .2in

\noindent \textbf{Theorem 2.} \textit{For each $\alpha \in \mathfrak{b}(M)$ there is a non-negative integer $\varepsilon (\alpha )$, which is determined by oriented intersection numbers in $M$, and an isomorphism of $R$-modules:
$$\mathcal{S}(M)\cong \bigoplus_{\alpha \in \mathfrak{b}(M)}R/(q^{2\varepsilon(\alpha)}-1).$$}

The isomorphism of the theorem is determined by a choice of standard links with given free homotopy classes of components. A precise description of the \textit{index} $\varepsilon (\alpha )$ of $\alpha $ is given below.

Note that each $0$-dimensional twisted homology module is a direct sum of cyclic modules with cyclic summands corresponding to path components [W]. The theorems above provide a description of the relations in the skein module through monodromies of paths in $I(M)$. The skein module $\mathcal{S}(M)$ is free if and only if all monodromies are trivial, and $\mathcal{S}(M)\cong SR\hat{\pi}(M)$ in this case.

The image of a link $K$ in the skein module is determined by the free homotopy classes of its components (some element $\alpha \in \mathfrak{b}(M)$ determining a path component of $\Lambda (M)$) and some integer number 
$w(K)$, which is well-defined modulo $2\varepsilon (\alpha )$. It follows from the skein relations that this number can be interpreted as the writhe (relative to a chosen standard link) of $K$. Thus, by computing the skein module, the indeterminancy of the universal writhe invariant is determined. But it should be noted that the writhe defined in this way does \textit{not} satisfy the usual behaviour with respect to \textit{smoothing} of crossings as known for $S^3$. This \textit{homological} writhe of links is measured in the quantum deformation of the first homology group of $M$ [P3].

\vskip .1in

Now we describe the index function $\varepsilon $. 
Throughout we will use the natural Hurewicz isomorphism from oriented singular bordism to oriented homology in dimensions $1$ and $2$ (see e.\ g.\ ([Ka], 13.15)). Let 
$$\lambda : H_2(M)\otimes H_1(M)\longrightarrow \mathbb{Z}$$
be defined by the oriented intersection number of oriented (maybe singular) closed surfaces and loops in $M$. 
Let
$$\mu :\mathfrak{b}(M)\rightarrow H_1(M)$$
be defined by taking the sum of the homology classes resulting by application of the Hurewicz homomorphism to the free homotopy classes in $\alpha \in \mathfrak{b}(M)$.
For $\beta \in \hat{\pi}(M)$ let $f_{\beta}: S^1\rightarrow M$ be a representing map. Let
$$\mu : \pi_1(\mathbb{L}(M),f_{\beta})\rightarrow H_2(M)$$  
be defined by taking the homology class of the map
$c' :S^1\times S^1\rightarrow M$ adjoint to $c\in \pi_1(\mathbb{L}(M),f_{\beta})$.
Now for $\alpha =[\alpha_1,\ldots ,\alpha_r]\in \mathfrak{b}(M)$ let $\Gamma (\alpha )$ be the subgroup of $\mathbb{Z}$ generated by all elements in 
$$\lambda (\mu(\pi_1(\mathbb{L}(M),f_{\alpha_i})),\mu(\alpha))$$
for all $1\leq i\leq r$.  If  $\Gamma (\alpha )$ is nontrivial let $\varepsilon (\alpha )$ be a positive generator, otherwise let $\varepsilon(\alpha )$ be $0$. Note that $\varepsilon (\emptyset )=0$, where $\emptyset $ is the unique sequence in $\mathfrak{b}(M)$ of length $0$. Obviously $\varepsilon (\alpha )$ does not depend on the choice of maps $f_{\beta }$.    

\vskip .1in

The following result can be considered as the singular version of Przytycki's theorem, \textbf{(a)} (also the proof gives the idea for \textbf{(b)}):

\vskip .1in

\noindent \textbf{Proposition.} \textit{$\mathcal{S}(M)$ is isomorphic to $SR\hat{\pi}(M)$ (equivalently is free) if and only if each mapping from a torus to $M$ is homologous into $\partial M$.} 

\vskip .1in

\noindent \textit{Proof.}  If each mapping from a torus to $M$ is homologous into $\partial M$ then all oriented intersection numbers of closed curves and tori are trivial. Thus the index is trivial for each $\alpha \in \mathfrak{b}(M)$. Assume that there is a mapping from a torus, which is not homologous into $\partial M$. It follows by Poincare duality (see ([K2], App.\ A)) that there exists an oriented loop, which has nontrivial intersection number with the singular torus. The map from the torus is the adjoint map of a representative of some element in $\pi_1(\mathbb{L}(M),f_{\beta})$ for some $\beta \in \hat{\pi}(M)$. Assume that the intersection number of the torus map and $f_{\beta}$ is non-trivial. Then each knot with homotopy class $\beta $ is torsion in $\mathcal{S}(M)$ by theorem 2. Otherwise we can find $\beta '\in \hat{\pi}(M)$ such that $\beta '$ has non-trivial intersection number with the torus map.
Then for $\alpha :=[\beta, \beta ']$ we have $\varepsilon(\alpha )\neq 0$. Thus $\mathcal{S}(M)$ has torsion. $\square$.

\section{The general results}

The quantum deformation of the fundamental group determines the writhe invariants of links in oriented $3$-manifolds in a universal way. The writhe of a link in $S^3$ naturally is the sum of the self-writhe and the total linking number (sum of pairwise linking numbers). Accordingly one can define a 
skein module of oriented $3$-manifolds, which contains the self-writhe and total linking number invariants. It turns out that the computation of this module does not require ideas beyond those used for the computation of the module $\mathcal{S}(M)$.

Let $R':=\mathbb{Z}[q_1^{\pm 1},q_2^{\pm 1}]$. Let $\mathcal{S}'(M)$ be the quotient of $R'\mathfrak{L}(M)$ by the submodule generated by all elements of the form $K_+=q_1^2K_-$ resp.\ $K_+=q_2^2K_-$ for crossings of same resp.\ distinct components, and $K^{(1)}=q_1K$. Using the coefficient homomorphism $R'\rightarrow R$, which maps $q_1$ and $q_2$ to $q$, Pryztycki's universal coefficient theorem [P2] induces the natural isomorphism
$$\mathcal{S}(M)\cong \mathcal{S}'(M)\otimes_{R'}R .$$
Similarly, the coefficient homomorphism $R'\rightarrow R$, which maps $q_1$
to $1$ and $q_2$ to $q$ induces the homomorphism onto the linking number skein module $\mathcal{L}(M)$ considered in [K1]. Finally there is the naturally homomorphism of skein modules defined from the ring homomorphism $R'\rightarrow R$, which maps $q_2$ to 1 and $q_1$ to $q$. Let $\mathcal{W}(M)$ denote the resulting skein module. In this skein module the \textit{relative self-writhe} of links is measured. We will state some detailed results about this module in theorem 5 and the following corollaries. The modules discussed in this paper can be summarized in the following diagram of natural $(R',R)$-epimorphisms:

\[
\begin{array}{ccccc}
\mathcal{L}(M) & \longleftarrow & \mathcal{S}'(M) & \longrightarrow & \mathcal{W}(M) \\
& & \Big\downarrow & & \\
& & \mathcal{S}(M) & & 
\end{array}
\] 

\noindent The modules $\mathcal{L}(M)$ resp.\ $\mathcal{W}(M)$ contain the 
\textit{link homotopy} resp.\ \textit{anti-link homotopy} part of $\mathcal{S}'(M)$ in a natural way.

\vskip .1in

The following two results are immediate generalizations of theorems 1 and 2. Let $\mathbb{N}$ denote the set of non-negative integers. 
 
\vskip .1in

\noindent \textbf{Theorem 3.} \textit{There exists a system of local coefficients $\mathfrak{Z}'$ in $I(M)$ with $\mathfrak{Z}'_x\cong R'$ such that
$$H_0(I(M),\mathfrak{Z}')\cong \mathcal{S}'(M).$$}

\noindent \textbf{Theorem 4.} \textit{For each $\alpha \in \mathfrak{b}(M)$ there exists an index
$$\varepsilon '(\alpha )=(\varepsilon_1(\alpha),\varepsilon_2(\alpha),\varepsilon_3(\alpha )\in \mathbb{N}\times \mathbb{Z}\times \mathbb{N},$$
determined by oriented intersection numbers as described below,
such that 
$$\mathcal{S}'(M)\cong \bigoplus_{\alpha \in \mathfrak{b}(M)}R'/(q_1^{2\varepsilon_1(\alpha)}q_2^{2\varepsilon_2(\alpha)}-1,q_1^{2\varepsilon_3(\alpha)}-1).$$}
 
For $\alpha =[\alpha_1,\ldots,\alpha_r]$ and $f_{\alpha_i}$ representing $\alpha_i$ let $\Gamma'(\alpha )$ be the subgroup of $\mathbb{Z}\times \mathbb{Z}$ generated by all elements (omitting the Hurewicz homomorphisms from the notation):
$$(\lambda(c,\alpha_i),\lambda(c,\alpha\setminus \alpha_i)),$$
where $c$ runs through $\pi_1(\mathbb{L}(M),f_{\alpha_i})$ and $1\leq i\leq r$.
Then $\Gamma '(\alpha )$ is generated by two elements $(\varepsilon_1(\alpha),\varepsilon_2(\alpha))$ and $(\varepsilon_3(\alpha),0)$ with $\varepsilon_1(\alpha),\varepsilon_3(\alpha)\in \mathbb{N}$. ($\Gamma '(\alpha )$ is cyclic if and only if
$\varepsilon_3(\alpha)=0$.)

\vskip .1in

It follows from the skein relations that the image of a given link $K$ in $\mathcal{S}'(M)$ can be written as $q_1^{w_1(K)}q_2^{w_2(K)}$ multiplied by a standard link with homotopy classes determined by some $\alpha \in \mathfrak{b}(M)$. Here $w_1(K)$ and $w_2(K)$ are the relative self-writhe and total linking number of $K$. The pair $(w_1(K),w_2(K))\in \mathbb{Z}^2$ is well-defined modulo the subgroup $\Gamma'(\alpha)$. Obviously 
the sum of the total linking number and the self-writhe is the writhe discussed before.
The number $\varepsilon_2(\alpha )$ is the linking number index as defined in [K1]. The $gcd$ of $(\varepsilon_1(\alpha )+\varepsilon_2(\alpha))$ and $\varepsilon_3(\alpha )$ is just $\varepsilon (\alpha )$ as defined in section 1. 

\vskip .1in

\noindent \textbf{Proposition.} \textit{$\mathcal{S'}(M)$ is isomorphic to $SR'\hat{\pi}(M)$ if and only if each mapping from a torus to $M$ is homologous into $\partial M$.}

\vskip .1in

\noindent \textit{Proof.} The only if part follows from Przytycki's universal coefficient result and the proposition in section 1. If every singular torus is homologous into $\partial M$ then all intersection numbers of singular tori and loops are zero. So the subgroups $\Gamma '(\alpha )\subset \mathbb{Z}\times \mathbb{Z}$ are trivial for all $\alpha $. It follows that $\mathcal{S}'(M)$ is free.
$\square$ 

\vskip .1in

In section 6 we will show that the relations in the natural presentation of the module $\mathcal{W}(M)$, equivalently the indeterminancies of the relative self-writhe invariants of links in $M$, are determined only from intersection numbers of \textit{$2$-spheres} and loops. More precisely the following result is shown:

\vskip .1in 

\noindent \textbf{Theorem 5.} \textit{For $\beta \in \hat{\pi}(M)$ let $\mu (\beta )$ be the smallest positive intersection number with a singular $2$-sphere, or $0$ if all intersection numbers are trivial. Then there is an isomorphism
$$\mathcal{W}(M)\cong \bigoplus_{\alpha\in \mathfrak{b}(M)}R/(q^{2\mu(\alpha)}-1),$$
where $\mu (\alpha )$ is the $gcd$ of the $\mu (\alpha _i)$ for $\alpha =[\alpha_1,\cdots,\alpha_r]\in \mathfrak{b}(M)$.}

\vskip .1in

\noindent \textbf{Remarks.} \textbf{(a)} Obviously $\mu (\alpha )$ only depends on the homology classes of the $\alpha_i\in \hat{\pi}(M)$ in $\alpha $.
  
\noindent \textbf{(b)} It follows from the sphere theorem (see e.\ g.\ [J], I.\ 9.\ ) that if a loop $\gamma $ is intersected by a singular $2$-sphere with non-trivial intersection number then there is also some \textit{embedded} $2$-sphere intersecting that loop with non-trivial intersection number. (Let $N_{\gamma}$ be the $\pi_1(M)$-invariant subgroup of $\pi_2(M)$ consisting of all those classes in $\pi_2(M)$, which have trivial intersection number with $\gamma $. Then given the assumption there exists an embedded $2$-sphere not in $N_{\gamma }$.)
  
\vskip .1in

\noindent \textbf{Corollary 1.} \textit{The module $\mathcal{W}(M)$ is isomorphic to $SR\hat{\pi}(M)$ if and only if there are no non-separating $2$-spheres in $M$ (equivalently each mapping of a $2$-sphere to $M$ is homologous into $\partial M$). }

\vskip .1in

If  $\mathcal{W}(M)$  is free then the relative writhe of \textit{knots} is defined without any indeterminancy for all free homotopy classes of loops in $M$.

\vskip.1in

\noindent \textbf{Corollary 2 (Chernov, Hoste and Przytycki).} \textit{The framing of an oriented knot can be changed by a self-isotopy of the underlying oriented link if and only if $M$ has a separating $2$-sphere.}

\vskip .1in

This result has been announced in ([HP], p.\ 488), and is proved in [C]. For the proof of corollary 2 we only have to remark that if there exists a non-separating 2-sphere then the construction in section 5 performed on a knot actually gives rise to an isotopy of knots.  

\vskip .1in

Our results show that freeness of the modules $\mathcal{S}'(M)$, $\mathcal{S}(M)$ or $\mathcal{L}(M)$ is equivalent to the geometric condition that each singular torus is homologous into $\partial M$. (For $\mathcal{L}(M)$ see ([K1], theorem 1.2)). The freeness of $\mathcal{W}(M)$ is equivalent to the condition that each singular $2$-sphere is homologous into $\partial M$. Note that the quantization of the first homology group of $M$ is free if and only if \textit{all} surfaces (and thus also all singular surfaces) are homologous into $\partial M$ [P3] (equivalently $M$ is a submanifold of a rational homology $3$-sphere ([K2], App.\ A)). It seems to be an interesting question whether it is possible to describe the geometric condition that all singular surfaces of genus $\leq g$
are homologous into $\partial M$ in terms of link theory in $M$.      

\section{Skein relations are monodromies}

In this section we prove theorem 3. Many of the arguments used here are similar to those in [K1].

\vskip .1in

Recall that elements of $I(M)$ are unordered oriented immersions of circles in $M$. For a given element $x\in I(M)$ consider all unordered framed immersions with the same underlying oriented immersion $x$. We will consider two such immersions equivalent if (i) their framings are homotopic or (ii) if they differ by twists in the framings of the components, such that the total signed number of these twists is zero. An equivalence class of framed immersions like this is called a total framing of $x$. We define $\mathfrak{Z}'_x$ resp.\ $\mathfrak{Z}_x$ to be the free $\mathbb{Z}[q_2^{\pm}]$-module resp.\ the free abelian group generated by all total framings of $x$. Note that there is a transitive action of $\mathbb{Z}$ on the set of total framings of $x$ for each $x\in I(M)$. 

\vskip .1in

There is an obvious notion of isotopy of totally framed links in $M$.
Let $\mathfrak{T}(M)$ denote the corresponding set of isotopy classes. It is immediate from the definitions that the skein module $\mathcal{S}'(M)$ can also be defined from the set $\mathfrak{T}(M)$ using the same relations as in the original definition of $\mathcal{S}'(M)$. Note that the total framing of the link $K^{(1)}$ is given by the action of $+1$ on the total framing of $K$. 

\vskip .1in

The collection of modules $\mathfrak{Z}'_x$ forms a bundle of modules in $I(M)$.
The corresponding local system of coefficients in $I(M)$ is defined by describing, for each $\alpha \in \mathfrak{b}(M)$, the action of $\pi_1(I(M),x)$
on  $\mathfrak{Z}'_x$. Here $x=x_{\alpha}$ is a link with free homotopy classes of components given by $\alpha \in \mathfrak{b}(M)$.

\vskip .1in

 By fixing an order of link components we define covering spaces $\widetilde{I(M)} \rightarrow I(M)$ and $\widetilde{\Lambda (M)}
\rightarrow \Lambda (M)$ in the obvious way. We will have to choose the base points $x=x_{\alpha }$ in a specific way. If $\alpha $ contains a conjugacy class repeatedly then we can choose $x_{\alpha}$ to be symmetric in the following way: All components with the same free homotopy class will be \textit{parallel} to each other with respect to some framing. It follows that there are isotopies (supported in a neighboorhood of the link), which change the order of those components. Note that by composition with loops of this form we can achieve that each loop in $I(M)$ lifts to a loop in $\widetilde{I(M)}$. The loops in $I(M)$ defined by those isotopies are called \textit{small}.
Note that the bundle of modules lifts to the covering space. 

\vskip .1in
  
Now for $\gamma \in \pi_1(I(M),x)$ let $\tilde \gamma$ be the loop, which results, possibly after composition with a small loop, by lifting to $\widetilde{I(M)}$. We choose a total framing on $x$ and lift this total framing to the corresponding basepoint on $\tilde \gamma $. It follows from Hirsch theory that the framing canonically transports along $\tilde \gamma $ and defines a new total framing on $x$. It is easily seen to be well-defined and independent of all choices made. The difference to the original framing is determined by some element $k_1(\gamma )\in \mathbb{Z}$. 

\vskip .1in

Next perturb the loop $\tilde \gamma $ to be transversal so that all singular maps along the path are immersions with a single double point. 
Let $k_2(\gamma )$ resp.\ $k_3(\gamma )$ be the signed sum of those immersions along $\gamma $ (or $\tilde \gamma $) (i.\ e.\ positive resp.\ negative crossing changes), where the double point of the immersion is a self-crossing resp.\ a crossing of distinct components.

\vskip .1in

\noindent \textbf{Theorem 6.} \textit{The collection of homomorphisms (of groups)
$$k=k_{\alpha}: \pi_1(I(M),x_{\alpha})\rightarrow Aut(R'),$$
defined by: $k(\gamma )$ is multiplication by
$$q_1^{k_1(\gamma )+2k_2(\gamma )}q_2^{2k_3(\gamma )},$$
defines a system of local coefficients $\mathfrak{Z}'$ in $I(M)$.}

\vskip .1in

\noindent \textit{Proof.} We have to show that the maps $k$ are constant under homotopies of loops in $I(M)$. The homomorphism property follows easily from the definitions. We will show more generally that $k$ is well-defined on 
$\pi_1(\Lambda (M),x)$. (This will be important later on). 
 We can assume that $\gamma $ lifts to a loop in $\widetilde{I(M)}$ and consider a loop $\gamma '$ homotopic to $\gamma $. The homotopy lifts to a homotopy in $\widetilde{I(M)}$ and we can apply Lin's transversality results [L] to change the homotopy to transversal position with respect to the complex of singular links in $M$. Then the preimage of the set of singular links is a $1$-dimensional complex in the domain of the homotopy (annulus), with vertices in the boundary and interior vertices of valence $4$ or $1$. Those of valence $4$ describe immersions with two double points, those of valence $1$ describe kink resolutions
(for more details see [KL] or [K1]). In order to show that the contributions from the two loops are equal deform one into the other one across the annulus in a finite number of steps, crossing only over a single interior vertex at a time. It is easy to see that $k$ does not change by crossing over a vertex of valence $4$ because the corresponding contribution around the vertex is zero. If we cross over a vertex of valence $1$ we encounter a crossing change, which does not change the oriented isotopy type but does change the total framing. The framing change cancels the singularity contribution. Thus crossing over any vertex does not change the value of $k$. $\square$    

\vskip .1in

By theorem 6 the homology module $H_0(I(M),\mathfrak{Z}')$ is defined. In theorem 7 we will define the isomorphism with $\mathcal{S}'(M)$. Recall that, as a twisted $0$-dimensional homology module, $\mathcal{S}'(M)$ is a direct sum of cyclic modules $R'/I(\alpha)$, where $\alpha $ runs through the set of path components of $I(M)$ (which is in 1-1 correspondence with $\mathfrak{b}(M)$). Moreover the ideal $I(\alpha )$ is generated by all elements of the form $k(\gamma )-1$, where $\gamma $ runs through all elements of $\pi_1(\Lambda (M),f_{\alpha})$
(compare e.\ g.\ ([W], theorem 3.2)).
Thus the proofs of theorem 4 resp.\ 2 reduce to the problem of relating the monodromies to the intersection numbers as explained in sections 1 and 2. This will be done in section 4.  

\vskip .1in

\noindent \textbf{Theorem 7.} \textit{There is a natural isomorphism
$$\mathcal{S}'(M)\longrightarrow H_0(I(M),\mathfrak{Z}').$$}
\noindent \textit{Proof.} There is the obvious map $\mathfrak{T}(M)\rightarrow H_0(I(M),\mathfrak{Z}')$ defined by mapping a link to the $0$-chain represented by a representative link $x$. Note that the framing determines an element in the fibre over $x$. So we have defined a $0$-cycle for homology with twisted coefficients, see ([W], p.\ 266). Recall that $H_0(I(M),\mathfrak{Z}')$ is an $R'$-module (this can be seen e.\ g.\ from the definition as a homology group of the chain complex $C_*(\overline{I(M)})\otimes R'$, where $\overline{I(M)}$ is the universal covering space of $I(M)$, and the tensor product is over the action of the fundamental group defining the system of local coefficients). The map above extends uniquely to an $R'$-homomorphism
$R'\mathfrak{T}(M)\rightarrow H_0(I(M),\mathfrak{Z}')$. This homomorphism is onto because every $0$-cycle can be represented, up to homology, by a link in $M$. It remains to prove that the skein relations are contained in the kernel of this homomorphism. It is shown in ([W], p.\ 263) how the action of the fundamental group determines a bundle of modules in $I(M)$, from which the boundary operator is defined. It follows from the definitions preceding theorem 6 that skein equivalent links are homologous in the twisted homology module. $\square$

\section{Computation of the monodromies}

Recall from the proof of theorem 6 that the images $k(\gamma )$ only depend on the homotopy class of $\gamma $ in $\Lambda (M)$. Moreover, by using the \textit{small} loops of embeddings, which change the order of components (described before theorem 6), we can assume that $\gamma $ actually lifts to a loop in $\widetilde{\Lambda (M)}$. But $\pi_1(\widetilde{\Lambda (M)}, (f_1,\ldots,f_r))$ is isomorphic to the direct product of groups $\pi_1(\mathbb{L}(M),f_i)$. Because of the homomorphism property of $k$ (see theorem 6) it suffices to compute $k(\gamma )$ on those loops in $I(M)$, which fix all but one component. 
For such a special loop in $I(M)$
let $\alpha _i$ be the free homotopy class of the \textit{nontrivial} component(which actually moves in $M$ along $\gamma )$. As before let $\alpha_i$ be represented by $f_{\alpha_i}: S^1\rightarrow M$ and let $c\in \pi_1(\mathbb{L}(M),f_{\alpha_i})$ be a loop in the free loop space based in $f_{\alpha }$ and consider a loop in $I(M)$ defined in this way.

Next consider the trace of the nontrivial component of such a loop, which is an immersion 
from $S^1\times I$ to $M \times I$.
The immersions over the two boundary components can be glued together to yield an immersion (after perturbation near the boundary) 
$$f: S^1\times S^1 \rightarrow M\times S^1,$$
which is framed except along the gluing circle. The difference of the framings 
on the top and bottom is precisely given by $k_1(\gamma )$ and is the normal Euler number $\chi (f)$ of the immersion. This is well-known from the singularity interpretation of the normal Euler number by using a section of the $2$-dimensional normal bundle, which is generic with respect to the zero-section. 
 
Then $k(\gamma )$ is obviously computed from the self-intersections of this immersion and from the intersection numbers of the projection to $M$ with the \textit{constant} components. The oriented intersection numbers of the torus in $M$ with the other components determine $k_3(\gamma )$. So we only have to consider the immersed torus. 

It is a classical result, in its original version due to Whitney [Wh] (discussed e.\ g.\ in the work of Lashof and Smale [LS]), that the homological intersection number of $[f]\in H_2(M\times S^1)$ satisfies the following relation:
$$[f]\cdot [f]=2D(f)+\chi (f).$$
 Here $D(f)$ is the oriented self-intersection number of the immersion $f$, which is equal to $k_2(\gamma )$ by definition. 

Thus the proof of theorem 4 is complete by proving the 

\vskip .1in

\noindent \textbf{Lemma.} \textit{The homological intersection number $[f]\cdot [f]$ is equal to $2\lambda (c,\alpha_i)$.} 

\vskip .1in

It follows that
$$2\lambda (c,\alpha_i)=[f]\cdot [f]=k_1(\gamma )+2k_2(\gamma ).$$
By considering corresponding powers of $q_1$ resp.\ $q_2$
 we see that $\varepsilon ' $ is determined as described in theorem 4. 

\vskip .1in

\noindent \textit{Proof of Lemma.} Decompose $[f]$ by the K\"{u}nneth theorem
$H_2(M\times S^1)\cong H_2(M)\oplus (H_1(M)\otimes H_1(M))$ into a sum $a+b$ with $a\in H_2(M)$ and $b\in H_1(M)\otimes H_1(M)$ (using obvious identifications). Note that $a\cdot a=b\cdot b=0$ (intersection numbers in $M\times S^1$). Thus $[f]\cdot [f]=2a\cdot b$. But if $a$ comes from the homology class of a singular torus in $M$ (by inclusion into $M\times S^1$) and $b$ is of the form $\mu(\alpha_i)\times 1$, where $1$ is the generator of $H_1(S^1)$, then $a\cdot b$ is the intersection number of $a\in H_2(M)$ and $\alpha_i$ in $M$. In order to see that $[f]-b$ is equal to $a$ we use an easy bordism argument. $\square$      

\section{Computational tools and Proof of Przytycki's result}

Let $\Omega M$ be the based loop space of $M$. Consider the fibration:
$$\Omega M\overset{i}{\rightarrow} \mathbb{L}M\overset{p}{\rightarrow}M,$$
where $i$ is the inclusion and $p$ is the evaluation at the basepoint (see [V2] and [K1]). For given $b\in \pi_1(M)$ let
let $f_b$ be a basepoint in $\Omega (M)$ corresponding to $b$ and let $\beta \in \hat{\pi}(M)$ be the corresponding free homotopy class. 

There is the exact sequence of homotopy groups:
$$\pi_1(\Omega M,f_b)\overset{i_*}{\longrightarrow} \pi_1(\mathbb{L}M,f_b)\overset{p_*}{\longrightarrow}\pi_1(M)\overset{[\  ,b]}{\longrightarrow}\pi_1(M).$$

\vskip .1in

Note that there are isomorphisms $\pi_1(\mathbb{L}(M),f_b)\cong \pi_1(\mathbb{L}M,f_{\beta})$, where the map $f_{\beta }: S^1\rightarrow M$ is the map chosen in section 1.

\vskip .1in

There is an obvious isomorphism $\pi_2(M)\cong \pi_1(\Omega M,f_{\alpha})$.
(Deform a loop in $f_b$ in $\Omega (M)$ such that all the mappings $S^1\rightarrow M$ along the loop are constant on a fixed neighbourhood of the basepoint. We can assume that $f_b$ is chosen to have this property. Then the difference to the \textit{constant loop in $f_b$} is uniquely determined by some element in $\pi_2(M)$.)

\vskip .1in

Now the image $i_*(c)$ for $c\in \pi_2(M)$ can be described as follows. Represent $c$ by a mapping of a $2$-sphere $g$, which is transverse to $f_b$. Now deform the restriction of $f_b$ on a neighbourhood of the basepoint in $*\in S^1$ as follows: Choose an arc joining $f_b(*)$ with $g(*)$. Then deform $f_b$ along this arc until a small neighbourhood of $*$ is mapped to $g(S^2)$. Now deform across the image of the $2$-sphere in the obvious way.
This is a description of $i_*(c)\in \pi_1(\mathbb{L}M,f_b)$. Note that if the intersection number of $f_b$ and $g$ is $k\in \mathbb{Z}$ then 
$$\lambda (\mu(i_*(c)),\mu (f_b))=k.$$
Note that this intersection number is trivial for all $2$-sphere mappings, which are homotopic to $\partial M$. 

\vskip .1in

\noindent \textbf{Example.} Let $M=S^2\times S^1$. Then $\hat{\pi}(M)=H_1(M)\cong \mathbb{Z}$ and $\mathfrak{b}(M)$ can be identified with the set of unordered sequences of integer numbers. Note that $p_*$ is onto and $\pi_1(\mathbb{L}M,f_b)$ is the union of sets $p_*^{-1}(b')$ with a transitive action of $\pi_2(M)$ on each of these sets. 

Since $\pi_1(M)$ is cyclic, a mapping from a torus (representing some $x\in \pi_1(\mathbb{L}M,f_b)$) can be homotoped to a map $f$ so that the restriction to the $1$-skeleton $S^1\vee S^1$ maps into the core $*\times S^1$. Now this restriction can be extended to a map of a torus into $*\times S^1$. Thus the collection of intersection numbers $\lambda(\mu (x),\mu (\alpha ))$, where $x$ runs through all of $p^{-1}(b')$, does not depend on $b'$. Here $\alpha \in \mathfrak{b}(M)$ is any sequence of free homotopy classes. So we can restrict the computation to $b'=0$ and compute all intersection numbers from singular tori resulting from elements in the image of $i_*$.  

Note that $\pi_2(M)$ is generated by multiples of the standard sphere $S^2\times *$. Moreover, a loop in $M$ has homology class $k\in \mathbb{Z}$ if and only if the loop intersects $S^2\times *$ with intersection number $k$.
It follows that for $\alpha \in \mathfrak{b}(M)$ 
$$\varepsilon (\alpha )=\varepsilon [\alpha_1,\ldots ,\alpha_r]=\sum_{i=1}^r\alpha_i.$$
The subgroup $\Gamma (\alpha )$ of $\mathbb{Z}\times \mathbb{Z}$ is generated by the elements $(\alpha_i, \sum_{j\neq i}\alpha_j)$. So if we consider $\alpha =[1,1,\ldots ,1]$ of length $r$, we have the contributions $(1,r-1)$ for all $i$. Then $\varepsilon '(\alpha )=(1,r-1,0)$ and the corresponding cyclic summand of the skein module is $R'/(q_1^2q_2^{2(r-1)}-1)$. For $\alpha =[1,2]$ the resulting subgroup of $\mathbb{Z}\times \mathbb{Z}$ is generated by $(1,2)$ and $(2,1)$. The summand of the skein module corresponding to this $\alpha $ is given by $R'/(q_1^4q_2^2-1,q_1^6-1)$. Even though the computation is easy for each given $\alpha $ there seems to be some tricky combinatorics involved in order to be able to state a general answer.   

\vskip .2in 

Now we give a proof of the theorem of J.\ Przyycki from theorem 2. 

\vskip .1in

First we show \textbf{(b)}. Note that $\lambda(\mu(i_*(c)),\mu(K))=k$ for $K$ a knot in $M$ intersecting the $2$-sphere with homotopy class $c\in \pi_2(M)$
with intersection number $k\in \mathbb{Z}$. This shows \textbf{(i)}. To prove \textbf{(ii)} we define the element $c\in \pi_1(\mathbb{L}M,f_{\beta})$, where $f_{\beta }: S^1\rightarrow M$ is a representative of the given non-contractible curve on the torus with free homotopy class $\beta \in \hat{\pi}(M)$, by a deformation of the non-contractible curve over the torus back into itself. Then let $\alpha = [\beta ,\gamma ]$ where $\gamma \in \hat{\pi}(M)$ is the free homotopy class of $K'$. It follows that $\lambda (\mu (c),\mu (\alpha ))=k$ because pushing the non-contractible curve away from the embedded torus shows that the intersection number of $\gamma $ with $c$ is zero. This proves \textbf{(ii)}. 

\vskip .1in
 
In order to prove part \textbf{(a)} of Przytycki's result we have to show that $\varepsilon (\alpha )=0$ if there are no non-separating $2$-spheres or tori in $M$. Conversely we claim that $\varepsilon (\alpha )\neq 0$ implies the existence of a non-separating $2$-sphere or torus.      
Now it is immediate from the definitions that there exists $\alpha \in \mathfrak{b}(M)$ satisfying $\varepsilon (\alpha )\neq 0$ if and only if there a singular torus and an oriented loop in $M$, which have non-trivial intersection number. Note that each singular $2$-sphere, by composition with the canonical projection $S^1\times S^1\rightarrow S^2$, induces a singular torus, preserving intersection numbers with a loop. Obviously (see e.\ g.\ ([K2], App.\ A)) each oriented \textit{embedded} surface with non-trivial intersection number with an oriented loop is non-separating.

Thus Przytycki's theorem follows from the following easy consequence of D.\ Gabai's fundamental results.

\vskip .1in

\noindent \textbf{Theorem.} \ \textit{Suppose that $M$ is a compact oriented $3$-manifold and let $\gamma $ be an oriented loop in $M$. If there is
 a singular torus in $M$ with non-trivial intersection number with $\gamma $ then there also exists an embedded $2$-sphere or torus with non-trivial intersection number with $\gamma $.}

\vskip .1in 

\noindent \textit{Proof.} We apply corollary 6.18.\ from [G] to the homology class $z\in H_2(M)$ determined by the singular torus. It follows that the singular Thurston norm of $z$ is zero, so by Gabai's result also the \textit{embedded} Thurston norm of $z$ is zero. It follows that $z$ can be represented by a disjoint union $T$ of embedded $2$-spheres and tori. This representative of $z$ has the same intersection number with $\gamma $ as the original singular torus. Obviously there is at least one component of $T$, which has non-trivial intersection number with $\gamma $. $\square$  

\vskip .1in

\noindent \textbf{Remark.} Note that it is not possible to conclude from singular to embedded tori and $2$-spheres \textit{preserving the intersection number}. In section 6 we will show that this is possible if we have the assumption of a singular $2$-sphere. 

\section{The indeterminancy of the self-writhe}

 Throughout we denote the homology class of a map $h$ from a sphere or torus to $M$ by $[h]$.

Consider a singular torus $f: T:=S^1\times S^1\rightarrow M$ with $f|\gamma =:g$ for some non-contractible curve $\gamma $ on $T$ such that $\lambda ([f],[g])=k\neq 0\in \mathbb{Z}$. We will show that there exists a singular $2$-sphere, which has intersection number $k$ with $g$. Then theorem 5 follows from theorem 4 and the results of section 5.    

First assume that $M$ is a connected sum of $3$-manifolds $M_1$ and $M_2$, or $M$ is $\partial$-compressible. In the second case $M$ is the result of attaching a $1$-handle $H$ (the neighbourhood of the attaching disk) along a pair of disks to a manifold $M'$. 
 We can assume that $f$ intersects the corresponding $2$-sphere resp.\ disk pair, denoted by $D$ in each case, transversely such that $f^{-1}(D)=:C$ is a disjoint union of closed curves on $T$. Suppose all components of $C$ are trivial on $T$. Then we can assume that $\gamma $ is disjoint from $C$. By the usual cut and paste arguments we construct a torus mapping to one of $M_1$ or $M_2$ resp.\ to $M'$, without changing the intersection number with $g$. (Note that the resulting mapping of torus cannot map into $H$ because intersection numbers of tori and curves in $H$ are trivial.) Suppose $C$ is a union of parallel copies of some non-contractible curve on $T$ (and possibly trivial curves). We do surgery on the torus mapping as before. Now the result is a mapping from a union of $2$-spheres to $M$. Note that in the process of glueing maps from disks (into $D$) along the boundaries of $T-N(C)$, where $N(C)$ is a neighbourhood of $C$ in $T$, the intersection number with $g$ is not changed. (In order to see this note that $f(C)$ is a collection of immersed circles in $D$ and we can assume that $g$ intersects $D$ only in a finite number of points away from $f(C)$. This shows that, when cutting the torus, the intersection number does not change. When disks are glued back in pairs with different orientations in the second step of the surgery the intersection number also remains unchanged.) Finally we can construct an \textit{oriented} connected sum of the resulting collection of maps of $2$-spheres and have found a mapping from a $2$-sphere to $M$ having intersection number $k$ with $g$. (Omit all components with trivial intersection number with $g$, the remaining $2$-spheres are non-separating so there is no problem in forming the oriented connected sum.)      

By iteration of the argument above it follows that  \textit{either} we find a $2$-sphere with homotopy class $c\in \pi_2(M)$ such that $\lambda (c,[g])=k$, or we can assume that $T$ maps into a $3$-manifold, which is irreducible and $\partial$-irreducible or is homeomorphic to $S^2\times S^1$.

In the last case the result follows from the considerations in the example in section 5.

So we can assume that $f$ maps into a $3$-manifold $M$, which is irreducible and $\partial$-incompressible. Since there is a non-trivial intersection number we know that $M$ is sufficently large. It follows from the main result of Jaco-Shalen-Johannson theory ([JS], p.\ 173) that $f$ is homotopic to a mapping into a Seifert fibred space $N$ contained in $M$. Moreover the boundary of $N$ is incompressible in $M$. Using this it follows easily from a theorem of Alexander ([R], p.\ 107) that $N$ is irreducible.
Now consider the closed oriented Seifert fibred space $N'$ defined by the double of $N$. It is not hard to show that $N'$ is irreducible. It follows from the parallizability of oriented $3$-manifolds and results of M.\ Hirsch that $f$ can be approximated by an immersion ([H], 5.7).
 Thus we can apply a result of P.\ Scott ([S], 6.4) to conclude that $f$ is homotopic, inside $N'$, to an immersion of a torus without triple points or to a composition of the projection $T$ onto the Klein bottle $K$ with some immersion $K\rightarrow M$. Then the result follows from the following lemma.

\vskip .1in

\noindent \textbf{Lemma.} \textit{Suppose that $f: T\rightarrow M$ is either an immersion without triple points or a composition of the double cover $T\rightarrow K$
with a mapping $K\rightarrow M$. Then $\lambda ([f],[g])=0$, where $g=f|\gamma $ and $\gamma $ is a non-contractible loop on $T$.} 

\vskip .1in

\noindent \textit{Proof.} In the second case the result follows easily from the different orientations in the two intersection points on $T$ corresponding to some intersection point of $K$ and $g$, and in fact holds for all loops in $M$.
 In the first case the self-intersection curves of $f$ consist of a disjoint union $C$ of non-contractible curves, all parallel to a single curve $\rho $, and a collection of trivial curves on $T$. Obviously $\lambda ([f],[\rho ])=0$. Let $\sigma $ be a non-contractible curve on $T$, which intersects $\rho $ in single point.
Then $\lambda ([f],[\sigma ])=0$ also holds. In fact, $f$ restricts to a \textit{trivial} double cover $C\rightarrow C'$. Moreover we can assume that $\sigma $ does not intersect any trivial components of $C$. Then the two contributions in $\lambda ([f],[\sigma ])$, which result from intersections of $\sigma $ with the two parallel copies of $\rho $ in $C$ ( mapping to the same circle in $C'$) cancel. This proves the result. $\square$

\vskip .3in

\noindent \textbf{Acknowledgements.} 

\vskip .1in

It is a pleasure to thank J.\ Przytycki for bringing the problem discussed in this paper to my attention, and for many stimulating discussions and generous hospitality.  

\vskip 0.3in

\begin{scriptsize}

\centerline{\textsc{REFERENCES}} 

\vskip .2in

\begin{description}

\item[[C]]
\textsc{V.\ Chernov}.
\textsl{Framed knots in $3$-manifolds}.
preprint math.GT/0105139, http://xxx.lanl.gov, 2001 

\item[[G]]
\textsc{D.\ Gabai}.
\textsl{Foliations and the topology of $3$-manifolds}.
J.\ Differential Geometry 18, 1983, 445--503

\item[[H]]
\textsc{M.\ W.\ Hirsch}.
\textsl{Immersions of manifolds}.
Transactions of the AMS 93, 1959, 242--276

\item[[HP]]
\textsc{J.\ Hoste, J.\ Przytycki}.
\textsl{Homotopy skein modules of orientable $3$-manifolds}.
Math.\ Proc.\ Camb.\ Phil.\ Soc.\ (11990), 108, 475--488

\item[[J]]
\textsc{W.\ Jaco}.
\textsl{Lectures on Three-Manifold Topology}.
Regional Conference Series in Mathematics Number 43, AMS 1980

\item[[JS]]
\textsc{W.\ H.\ Jaco, P.\ B.\ Shalen}.
\textsl{Seifert Fibred Spaces in $3$-manifolds}
Memoirs of the AMS 220, Vol 21, 1979

\item[[K1]]
\textsc{U.\ Kaiser}.
\textsl{Presentations of $q$-homotopy skein modules of oriented $3$-manifolds}.
Journal of Knot Theory and its Ramifications, no. 3, 2001, 461--491

\item[[K2]]
\textsc{U.\ Kaiser}.
\textsl{Link theory in manifolds}. 
Lecture Notes in Mathematics 1669, Springer Verlag 1997,
59--77

\item[[Ka]]
\textsc{L.\ Kauffman}.
\textsl{On Knots}.
Annals of Mathematical Studies 115, Princeton University Press 1987

\item[[KL]]
\textsc{E.\ Kalfagianni, X.\ S.\ Lin}.
\textsl{The HOMFLY polynomial for links in rational homology spheres}.
Topology 38, no. 1, 1999, 95--115

\item[[Ki]]
\textsc{R.\ Kirby}.
\textsl{Problems in low dimensional topology}.
in 1993 Georgia International Topology Conference Proceedings, Part II, ed.\ by W.\ Kazez,
AMS/IP Studies in Advanced
Mathematics 1997

\item[[LS]]
\textsc{R.\ K.\ Lashof, S.\ Smale}.
\textsl{Self-intersections of immersed manifolds}.
Journal of Mathematics and Mechanics, Vol.\ 8, N0.\ 1, 1959, 143--157

\item[[L]] 
\textsc{X.\ S.\ Lin}.
\textsl{Finite type link invariants of
$3$-manifolds}.
Topology 33, no.1, 1994, 45--71

\item[[P1]]
\textsc{J.\ Przytycki}.
\textsl{Algebraic topology based on knots: an introduction}.
in Proceedings of Knots 96, edited by Shin'ichi Suzuki, 1997, 279--297

\item[[P2]]
\textsc{J.\ Przytycki}.
\textsl{Skein modules of $3$-manifolds}.
Bull.\ Ac.\ Pol.\ Math.\ 39(1-2), 1991, 91-100

\item[[P3]]
\textsc{J.\ Przytycki}.
\textsl{A $q$-analogue of the first homology group of a $3$-manifold}.
Contemporary Mathematics Volume 214, 1998, 135--143

\item[[R]]
\textsc{D.\ Rolfsen}.
\textsl{Knots and Links}
Mathematics Lecture Series 7, Publish or Perish Inc.\ , Second Printing 1990

\item[[S]]
\textsc{P.\ Scott}.
\textsl{There are no fake Seifert fibre spaces with infinite $\pi_1$}.
Annals of Mathematics, 117 (1983), 35--70

\item[[V1]]
\textsc{V.\ Vassiliev}.
\textsl{On invariants and homology of spaces of knots in arbitrary manifolds}.
American Math.\ Soc.\ Translations (2), Vol.\ 185, 1998

\item[[V2]]
\textsc{V.\ Vassiliev}.
\textsl{Complements of discriminants of smooth maps: topology and applications}.
Transl.\ Math.\ Monographs, vol.\ 98, AMS, Providence, RI 1994

\item[[W]]
\textsc{G.\ W.\ Whitehead}.
\textsl{Elements of Homotopy Theory}.
Graduate Texts in Mathematics Vol.\ 61, Springer Verlag 1978

\item[[Wh]]
\textsc{H.\ Whitney}.
\textsl{The self-intersections of a smooth $n$-manifold in $2n$-space}.
Ann.\ Math.\ (1944), 220--240

\end{description}

\end{scriptsize}
\end{document}